\documentclass[12pt,a4paper]{amsart}
\usepackage{amssymb,latexsym,amsmath,amscd}
\newcommand{\Pn}{\mathbb{P}^{n}}
\newcommand{\PP}{\mathbb{P}}
\newcommand{\ZZ}{\mathbb{Z}}

\newcommand{\OO}{\mathcal{O}}

\newcommand{\NN}{\mathbb{N}}

\newcommand{\FF}{\mathcal{F}}

\newcommand{\JJ}{\mathcal{J}}

\newcommand{\lra}{\longrightarrow}
\newcommand{\ra}{\rightarrow}

\newcommand{\al}{\alpha}
\newcommand{\be}{\beta}
\newcommand{\Hom}{\textrm{Hom}}

\newcommand{\srel}{\stackrel}

\newcommand{\depth}{\operatorname{depth}}

\theoremstyle{definition}
\newtheorem{defi}{Definition}[section]
\theoremstyle{plain}
\newtheorem{lema}[defi]{Lemma}
\newtheorem{teo}[defi]{Theorem}

\newtheorem{prop}[defi]{Proposition}

\newtheorem{rk}[defi]{Remark}

\newtheorem{ex}[defi]{Example}

\newtheorem{nota}[defi]{Notation}
\newtheorem{teo-def}[defi]{Theorem/Definition}
\newenvironment{proofteo}{\noindent {\textit{Proof of Theorem \ref{MRC}.}}}{\quad \hfill $\Box$}
\newenvironment{prooflema1}{\noindent {\textit{Proof of Proposition \ref{firstlink}.}}}{\quad \hfill $\Box$}
\newenvironment{prooflema2}{\noindent {\textit{Proof of Proposition \ref{secondlink}.}}}{\quad \hfill $\Box$}
\newenvironment{prooflema3}{\noindent {\textit{Proof of Proposition \ref{thirdlink}.}}}{\quad \hfill $\Box$}
\newenvironment{prooflema4}{\noindent {\textit{Proof of Proposition \ref{fourthlink}.}}}{\quad \hfill $\Box$}

\title[The minimal resolution conjecture \dots]{The minimal resolution conjecture for points on the cubic surface}

\author{M. Casanellas}
\thanks{Partially supported by Ministerio de Educaci\'on y Ciencia of Spain, Programa Ramon y Cajal and BFM2003-06001.\\ To appear in Canadian Journal of Mathematics.}
\address{Departament Matem\`atica Aplicada  I \\ ETSEIB UPC. \\ Av. Diagonal 647 \\ 08028-Barcelona. Spain \\ marta.casanellas@upc.edu}
\begin{document}
\maketitle
\begin{abstract}
In this paper we prove that the generalized version of the Minimal Resolution Conjecture stated in \cite{Farkas} holds for certain general sets of points on a smooth cubic surface $X \subset \PP^3$. The main tool used is Gorenstein liaison theory and, more precisely, the relationship between the free resolutions of two linked schemes.
\end{abstract}
\section{Introduction}
The Minimal Resolution Conjecture for points in projective spaces was first stated by Lorenzini in \cite{Lorenzini}. Roughly speaking it says that the graded minimal free resolution of a general set of $t$ points in $\PP^n$ has no \emph{ghost} terms. This conjecture is known to hold for $n \leq4$ (\cite{Gae}, \cite{BG} and \cite{Walter}) and for large values of $t$ for any $n$ (cf. \cite{HS}) but it does not hold in general: Eisenbud, Popescu, Schreyer and Walter proved that it fails for every $n \geq 6$, $n \neq 9$ (see \cite{EPSW}). Farkas, Musta\c{t}\v{a} and Popa introduced a generalized version of this conjecture for points in arbitrary projective varieties (see \cite{Mustata} and \cite{Farkas}). Namely, if $X \subset \PP^n$ is any projective variety, the Minimal Resolution Conjecture for a general set of $t$ points $Z$ on $X$ predicts that the Betti numbers of the ideal of $Z$ are completely determined by those of the ideal of $X$ (see \cite{Farkas} and section 2 for a precise statement of the conjecture). When $X=\PP^n$, this formulation of the conjecture coincides with the conjecture formulated by Lorenzini.

This generalization of the Minimal Resolution Conjecture has been studied in \cite{Farkas} and has been proven to hold for a general set of points of any sufficiently large degree when $X \subset \PP^n$ is a canonical curve. However, it always fails for sets of points on curves of large degree (see \cite{Farkas}). Guiffrida, Maggioni and Ragusa (see \cite{Giuffrida}) proved that this conjecture also holds for any general set of points when $X$ is a smooth quadric surface in $\PP^3$. In this paper we study the Minimal Resolution Conjecture for general sets of points on a smooth cubic surface in $\PP^3.$ 

The goal of this paper is to prove that $t$ general points on a smooth cubic surface $X \subset \PP^3$ satisfy the Minimal Resolution Conjecture when $t$ is equal to $\frac{3}{2}a(a-1)+a$, $\frac{3}{2}a(a-1)+2a$, $\frac{3}{2}a(a-1)+a+1$ or $\frac{3}{2}a(a-1)+2a+1$ for some $a \in \ZZ$ (see Theorem \ref{MRC}). We also give the precise minimal free resolutions of these sets of points. All these sets of points happen to be \emph{level} so, in particular, this result provides proofs to some unjustified examples in \cite{GHMS} Appendix C (namely that 35,36,39,41 and 42 general points on the cubic surface are level). It also provides a partial answer to Question 7.5 of \cite{GHMS} ``\emph{If $X$ is an integral arithmetically Cohen-Macaulay variety of dimension $\leq 5$, and $Z$ is a general set of points on $X$, then does $Z$ have the expected resolution? In particular, is it true that $Z$ is level if it is numerically possible?}" and gives hope for an affirmative answer to this question.

Whereas \cite{Mustata} and \cite{Farkas} used Koszul cohomology to study the Minimal Resolution Conjecture for points on curves, the main tool used in this paper is Gorenstein liaison theory. R. Hartshorne proved in \cite{Hargor} that a general set of points in the cubic surface $X$ can be Gorenstein linked in $X$ to one point by a finite sequence of links. We make use of the $G$-links he describes to track the minimal free resolutions of the sets of points we are considering. It is a well known result in liaison theory that there is a relationship between the free resolutions of two linked schemes, but one needs to be careful because the minimality of the resolutions is not preserved in general.

Our initial motivation for studying the resolution of sets of points on the cubic surface was the connection between level sets of points and Ulrich bundles in the sense of \cite{CDH}, but this will be explained in a forthcoming paper.

The paper is organized as follows. In the next section we write down the Minimal Resolution Conjecture and we provide the results of Gorenstein liaison theory needed for the sequel. Section 3 is devoted to proving the main result of the paper.

\underline{Acknowledgments:} The author would like to thank R. Hartshorne for enjoyable discussions on sets of points on the cubic surface and thank him and J. Migliore for useful comments on a preliminary version of the present paper.

\section{Preliminaries}

Throughout this section $R$ denotes the polynomial ring in $n+1$
variables over an algebraically closed field $k$. Throughout all
the paper, for a subscheme $Z \subset \PP^n_k$ we denote by $I_Z
\subset R$ its saturated ideal, i.e. $I_Z=\oplus_{t \in \ZZ
}H^0(\PP^n, \JJ_Z(t))$. For any coherent sheaf $\FF$ over a
projective scheme $X \subset \PP^n$, $H^i_{\ast}(\FF)$ denotes the
sum $\oplus_{t \in \ZZ }H^i(X, \FF(t)).$ The regularity $\textrm{reg}(Z)$ of $Z$ is defined to be the regularity of $I_Z$, if $Z \neq \PP^n$ and 1 otherwise. The Hilbert polynomial of $Z$ will be denoted as $P_Z$.

The results of this paper deal with minimal free resolutions of
graded $R$-modules. We will always talk about \emph{graded}
minimal free resolutions over $R.$

\begin{defi} If $Z \subset \PP^n$ is a subscheme and
$$ d_{\bullet} \quad 0 \ra F_{n+1} \srel{d_{n+1}}{\ra} F_{n}
\srel{d_{n}}{\ra} \dots \ra F_1 \srel{d_{1}}{\ra} R
\srel{d_0}{\ra} R/I_Z \ra 0$$ is the minimal free resolution of
$R/I_Z$, the \emph{Betti numbers} $b_{i,j}(Z)$ are defined as
$$F_i=\oplus_{j \in \ZZ}R(-i-j)^{b_{i,j}(Z)}.$$
\noindent
The \textit{Betti diagram }of $Z$ has in the $(j,i)$-th position the Betti number $b_{i,j}(Z)$. The last nontrivial row of the Betti diagram of $Z$ is indexed by $\rm{reg}(Z)-1$ (see \cite{E}).

\end{defi}

\noindent \textbf{Minimality criterion:} We recall that a free
resolution
$$ d_{\bullet} \quad 0 \ra F_{n+1}
\srel{d_{n+1}}{\ra} F_{n} \srel{d_{n}}{\ra} \dots \srel{d_1}{\ra}
F_0 \srel{d_{0}}{\ra} M \ra 0$$ is minimal if, after choosing
basis of the modules $F_i$, the matrices representing the maps
$d_i$ do not have any non-zero scalar entry.

\noindent \textbf{Minimal Resolution Conjecture for points on
embedded varieties.}(see \cite{Farkas}) Let $X \subset \PP^n$ be
an irreducible projective variety of dimension $\geq 1$ and denote
by $P_X$ its Hilbert polynomial. Let $Z$ be a \emph{general} set
of $z$ points in $X$. In this paper we say that a statement holds
for a \emph{general} $Z$ if there is a nonempty open subset of the
family of all $Z$'s for which the statement holds.
We choose $r$ such that $P_X(r-1) \leq z < P_X(r)$ and assume that $r \geq \textrm{reg}(X)+1$. Then the
\emph{Minimal Resolution Conjecture} formulated in \cite{Farkas}
holds for the value  $z$ if for every set $Z$ of $z$ general
points on $X$,
$$b_{i+1,r-1}(Z)b_{i,r}(Z)=0 \textrm{ for all } i .$$

Musta\c{t}\v{a} \cite{Mustata} proved that the first rows of the Betti diagram of a general set of points $Z$ on $X$ coincide with the Betti diagram of $X$ and that there are two extra nontrivial rows at the bottom. He also gives lower bounds for the Betti numbers in these last two rows and the Minimal Resolution Conjecture states that these lower bounds are attained for a general set of points. More precisely, his results can be summarized as follows:
\begin{teo}\rm{(\cite[Theorem 1.2]{Farkas})}\label{general_results}
Assume that $Z \subseteq X$ is a general set of $z$ points,
with $P_X(r-1)\leq z<P_X(r)$ for some $r\geq m+1$, where $m={\rm reg}\,X$.
\item{\rm (i)} For every $i$ and $j\leq r-2$, we have $b_{i,j}(Z)
=b_{i,j}(X)$.
\item{\rm (ii)} $b_{i,j}(Z)=0$, for $j\geq r+1$ and there is an  
$i$ such that $b_{i,r-1}(Z)\neq 0$.
\item{\rm (iii)} If $d=\dim\,X$, then for every $i\geq 0$, we have
$b_{i+1,r-1}(Z)-b_{i,r}(Z)=Q_{i,r}(z)$, where
$$Q_{i,r}(z)=\sum_{l=0}^{d-1}
(-1)^l{{n-l-1}\choose{i-l}}\Delta^{l+1}P_X(r+l)-
{n\choose i}(z-P_X(r-1)).$$ In particular,
$b_{i+1,r-1}(Z)\geq {\rm max}\,\{Q_{i,r}(z),0\}$
and $b_{i,r}(Z)\geq {\rm max}\,\{-Q_{i,r}(z),0\}$.
\end{teo}

In other words, for $z \geq P_X(r-1)$ and $r \geq m+1=\rm{reg}(X)+1$, the Minimal Resolution Conjecture is satisfied if and only if
the minimal free resolution of $Z$ has no \textit{ghost} terms (i.e. there are no identical free summands in two consecutive steps of the resolution). Note that the Betti numbers of $X$ do not overlap with those of $Z$ because we are requiring $r \geq m+1$. Indeed, the last row in the Betti diagram of $X$ is indexed by $m-1$ while the part of the Betti diagram of $Z$ that is not determined by $b_{i,j}(X)$ occurs in the rows $r-1$ and $r$. Here is an example of how the Betti numbers of $X$ determine those of $Z$ if the Minimal Resolution Conjecture holds.

\begin{ex}\rm Let $Z$ be a set of $z=22$ general points on a smooth cubic surface $X \subset \PP^3$. As $P_X(t)= \frac{3}{2}t(t+1)+1$  and $\rm{reg}(X)=3$, we take $r=4$ so that $P_X(3)=19 \leq 22 <P_X(4)=31$. According to Theorem \ref{general_results}, the first 2 rows of the Betti diagram of $Z$ coincide with the Betti diagram of $X$ and in rows $3$ and $4$ we have $b_{1,3}-b_{0,4}=9$, $b_{2,3}-b_{1,4}=12$, $b_{3,3}-b_{2,4}=9$, $b_{4,3}-b_{3,4}=-3$. If the Minimal Resolution Conjecture holds for $Z$, one of the terms in each of these differences is 0. Therefore, as $b_{i,j}(Z)$ is always positive and $b_{i,j}(Z)=0$ for $i \geq 4$, the Betti diagram should be as follows
\begin{center}
\renewcommand{\arraystretch}{1.25}
\begin{tabular}[l]{c|c}
 & 0 1 2 3 \\
\hline
0 & 1 -- -- -- \\
1 & -- -- -- -- \\
2 & -- 1 -- -- \\
3 & -- 9 12 --  \\
4 & -- -- -- 3 
\end{tabular}
\end{center}
or, equivalently, the minimal free resolution should be
$$0 \ra R(-7)^{3} \ra R(-5)^{12}  \ra R(-4)^{9} \oplus R(-3)  \ra I_{Z}
 \ra 0 .$$
In Theorem \ref{MRC} we will prove that, indeed, 22 general points on a smooth cubic surface satisfy the Minimal Resolution Conjecture.
\end{ex}

\begin{ex}\label{exMustata}\rm Musta\c{t}\v{a} proved in \cite{Mustata} that when $z=P_X(r-1)$ or $z=P_X(r)-1$, the Minimal Resolution Conjecture holds for every $X$.
\end{ex}

When $X=\PP^n$ the Minimal Resolution Conjecture has been deeply
studied and is known to hold for $n \leq 4$ (see \cite{Gae},
\cite{BG} and \cite{Walter}). However, it fails in the middle of the resolution for every $n \geq
6$, $n \neq 9$. The conjecture is also known to hold for any set of points on a smooth quadric surface in $\PP^3$ (cf. \cite{Giuffrida}, although this paper preceeds \cite{Mustata}, they precisely prove that the minimal free resolution of any set of general points on a smooth quadric $X$ is completely determined by the Betti numbers of $X$). It is also known to hold, among other cases, for any set of general points on  a canonical curve $X$  and to fail when $X \subset \PP^n$ is a curve of large degree (cf. \cite{Farkas}).

The main tool we shall use to study the Minimal Resolution
Conjecture on the cubic surface is Gorenstein liaison theory.
A good reference for liaison theory is the book \cite{M}. Here we
just recall the definitions we need for this paper and some tools
related to Gorenstein liaison.

We recall that a scheme $V \subset \PP^N$ is an \textit{arithmetically
Cohen-Macaulay} scheme if its homogeneous coordinate ring is a
Cohen-Macaulay ring (i.e. $\dim{R/I_V}$ $= \depth{R/I_V}$). If $V$ is
of dimension $d \geq 1 $, $V$ is arithmetically Cohen-Macaulay if
and only if $H^i_{\ast}(\PP^n,\JJ_V)=0$ for $i=1, \dots, d.$ Any
0-dimensional scheme is arithmetically Cohen-Macaulay. An obvious
fact that we will use throughout the paper is the following: if
$C \subset \PP^n$ is arithmetically Cohen-Macaulay of dimension
$\geq 1$ and $Z \subset C$ is any subscheme, then
$I_{Z,C}=H^0_{\ast}(\JJ_{Z,C})=I_Z/I_C$.

A closed subscheme $V \subset \PP^n$ of codimension $c$ is {\em
arithmetically Gorenstein} if its coordinate ring $R/I_V$ is a Gorenstein $^*$local graded ring (in the sense of \cite{BroadmannSharp}). This is equivalent to saying that its saturated homogeneous ideal
$I_V$ has a graded minimal free $R$-resolution of the following
type: $$0 \ra R(-t) \ra F_{c-1} \ra \ldots \ra F_1 \ra F_0 \ra I_V
\ra 0.$$ In other words, $V \subset \PP^n$ is arithmetically
Gorenstein if and only if $V$ is arithmetically Cohen-Macaulay and
the last module in the minimal free resolution of its saturated
ideal has rank one.

\begin{defi}(see \cite{M})
We say that two subschemes $V_1$ and $V_2$ of $\Pn$ are
\emph{directly Gorenstein linked}, or simply \emph{directly
G-linked}, by an arithmetically Gorenstein scheme $G \subset \Pn$
if $I_G \subset I_{V_1} \cap I_{V_2}$ and we have $I_G: I_{V_1}=
I_{V_2}$ and $I_G: I_{V_2}= I_{V_1}.$ We say that $V_2$ is the
\emph{residual} to $V_1$ in $G$. When $G$ is a complete
intersection we talk about a $CI$-\textit{link}.
\end{defi}

Whenever $V_1$ and $V_2$ do not share any common component, the
fact of being directly $G$-linked by a scheme $G$ is equivalent to
saying that $G=V_1 \cup V_2.$ \emph{Gorenstein liaison} studies the equivalence relation generated by $G$-links.

One of the most interesting properties preserved by a $G$-link is
the free resolution. To explain how it works we first recall the mapping cone procedure.
\begin{lema}[Mapping cone procedure, cf. \cite{McLane}]\label{cone} Given a short exact sequence of finitely generated
$R$-modules
$$0 \lra A \srel{\al}{\lra} B \lra C \lra 0,$$
and free resolutions
$$ e_{\bullet} \quad 0 \ra G_{n+1} \srel{e_{n+1}}{\ra} G_{n} \srel{e_{n}}{\ra} \dots \srel{e_1}{\ra} G_0 \srel{e_{0}}{\ra} A \ra 0 $$
and
$$ d_{\bullet} \quad 0 \ra F_{n+1}
\srel{d_{n+1}}{\ra} F_{n} \srel{d_{n}}{\ra} \dots \srel{d_1}{\ra}
F_0 \srel{d_{0}}{\ra} B \ra 0,$$ then the map $\al$ lifts to a map
between the resolutions $\xi_{\bullet}: e_{\bullet} \lra
d_{\bullet}$ and a free resolution for $C$ is
$$0 \ra G_{n+1} \srel{c_{n+2}}{\ra} G_{n} \oplus F_{n+1} \srel{c_{n+1}}{\ra} \dots \srel{c_3}{\ra} G_{1} \oplus F_2 \srel{c_2}{\ra} G_0 \oplus F_1 \srel{c_1}{\ra} F_0 \srel{c_0}{\ra} C \ra 0$$
where $$c_{i+1}=\begin{pmatrix} -e_i & 0 \\
\xi_i & d_{i+1}
 \end{pmatrix} , \quad 1 \leq i \leq n.$$
\end{lema}

The resolution of $C$ produced in the Lemma \ref{cone} above is
not necessarily minimal even if those of $A$ and $B$ are.

In the following Lemma we recall how to pass from the free resolution of a scheme to the free
resolution of its residual in  an arithmetically Gorenstein
scheme.
\begin{lema}[see \cite{M}]\label{mapcone}
Let $Z, Z' \subset \PP^n$ be two arithmetically Cohen-Macaulay
subschemes of codimension $c$ directly G-linked by an
arithmetically Gorenstein scheme $G$. Let the minimal free
resolutions of $I_Z$ and $I_G$ be
$$0 \ra F_c \srel{d_c}{\ra} F_{c-1}
\srel{d_{c-1}}{\ra} \dots \ra F_1 \srel{d_{1}}{\ra} I_Z \ra 0$$
 and
$$ 0 \ra R(-t) \srel{e_c}{\ra} G_{c-1} \srel{e_{c-1}}{\ra} \dots \ra G_1 \srel{e_{1}}{\ra} I_G \ra 0 $$
respectively. Then
 the functor $Hom_R(\cdot,R(-t))$ applied to a free resolution of $I_Z/I_G$ gives a free resolution of
$Z'$. In particular,
$$ 0 \ra F_1^{\vee}(-t) \ra F_2^{\vee}(-t) \oplus G_1^{\vee}(-t) \ra \dots \ra F_c^{\vee}(-t) \oplus G_{c-1}^{\vee}(-t) \ra I_{Z'} \ra 0$$
is a free resolution of $I_{Z'}$ (not necessarily minimal)
obtained by mapping cone from the resolutions of $I_G$ and $I_Z$.
\end{lema}

\section{Minimal Resolution Conjecture on the cubic surface}
Let $X$ be a smooth cubic surface in $\PP^3$ defined by a homogeneous cubic form
$f$. From now on $R$ denotes the polynomial ring $k[x_0, \dots,
x_3]$ over an algebraically closed field $k$ and $R_X$ denotes the
ring $R_X=R/(f).$ When $Z$ is a subscheme of $X$, $I_Z$ denotes
its saturated ideal in $R$ and $I_{Z,X} \subset R_X $ the ideal
$I_Z/(f)=H^0_{\ast}(\JJ_{Z,X}).$

As the Hilbert polynomial of $X$ is $P_X(t)= \frac{3}{2}t(t+1)+1,$  $P_X(t)$ general points on
$X$ satisfy the Minimal Resolution Conjecture and so do $P_X(t)-1$
general points on $X$ (see Example \ref{exMustata}). We are going to prove that other families of
general points in $X$ satisfy the Minimal Resolution Conjecture.

\begin{nota}\rm We set the following notation. \begin{enumerate}
\item[(i)] If $t$ is a positive integer, $Z_t$ denotes a set of $t$ general points in $X$.
\item[(ii)] $m(a)=\frac{3}{2}a(a-1)+a$, $n(a)=\frac{3}{2}a(a-1)+2a$, $o(a)=\frac{3}{2}a(a-1)+a+1$, $p(a)=\frac{3}{2}a(a-1)+2a+1$ for any $a \in \ZZ$.
\item[(iii)] $C_0$ is any smooth conic on $X$, $\Gamma$ any twisted cubic on $X$ and
$L$ any of the 27 lines in $X$.
\item[(iv)] $H$ denotes a general hyperplane section of $X$.
\item[(v)] if $C$ is a curve on $X$, $H_C$ will be a general hyperplane section of $C$ and $K_C$ a canonical divisor on $C$.
\end{enumerate}
\end{nota}

The curves $C_0$, $\Gamma$ and $L$ are arithmetically
Cohen-Macaulay curves and, according to the proof of \cite{Hargor}
Proposition 2.4, any arithmetically Cohen-Macaulay curve on $X$ is
linearly equivalent to $C_0+aH_X$, $\Gamma+aH$, $L+aH$ or
$aH$, for some $a \in \NN$. The modules
$H^0_{\ast}(\PP^3,\JJ_{C_0,X})$,
$H^0_{\ast}(\PP^3,\JJ_{\Gamma,X})$ and
$H^0_{\ast}(\PP^3,\JJ_{L,X})$ are maximal Cohen-Macaulay
$R_X$-modules and their minimal free $R$-resolutions are:
$$0 \lra R(-3)^2 \srel{\varphi_1}{\lra} R(-1)\oplus R(-2) \lra I_{C_0,X} \lra 0,$$
$$0 \lra R(-3)^3 \srel{\varphi_2}{\lra}  R(-2)^3 \lra I_{\Gamma,X} \lra 0,$$
$$0 \lra R(-2) \oplus R(-3) \srel{\varphi_3}{\lra} R(-1)^2 \lra I_{L,X} \lra 0.$$
For certain matrices $\psi_j$, the pairs $(\varphi_j, \psi_j)$ are
matrix factorizations of $f$ (see \cite{Eisenbud}): $\varphi_j
\cdot \psi_j=f \cdot Id$, $\psi_j \varphi_j =f \cdot Id$,
$f=det(\varphi_j)=det(\psi_j)$, $j=1,2,3$.

The degree and genus of smooth arithmetically Cohen-Macaulay curves $C$ on $X$ are:
\begin{enumerate}
\item $d=3a-2$, $g=\frac{1}{2}(3a^2-7a+4)$ if $C \sim L+(a-1)H$.
\item $d=3a-1$, $g=\frac{1}{2}(3a^2-5a+2)$ if $C \sim C_0+(a-1)H$.
\item $d=3a$, $g=\frac{1}{2}(3a^2-3a)$ if $C \sim \Gamma+(a-1)H$.
\item $d=3a$, $g=\frac{1}{2}(3a^2-3a+2)$ if $C \sim aH$.
\end{enumerate}
En each case, one of these curves moves on a linear system of dimension $d+g-1$. Note that the degree and genus of an arithmetically Cohen-Macaulay curve $C$ determine the linear system in which $C$ belongs. Indeed, the four possibilities above imply that if two arithmetically Cohen-Macaulay curves on $X$ have the same degree and genus, they belong to the same linear system. 

The goal of this paper is to prove the following result.
\begin{teo}\label{MRC} Let $X \subset  \PP^3$ be a smooth cubic surface and consider the above notation.  Then, for $a \geq 3$, the graded minimal free resolutions of $m(a)$, $n(a)$, $o(a)$ and  $p(a)$ general points on $X$ are:
$$0 \ra R(-a-3)^{a-1} \ra R(-a-1)^{3a}  \ra R(-a)^{2a+1} \oplus R(-3)  \ra I_{Z_{m(a)}}
 \ra 0 ,$$
$$0 \ra R(-a-3)^{2a-1} \ra R(-a-2)^{3a}  \ra R(-a)^{a+1} \oplus R(-3)  \ra I_{Z_{n(a)}}
 \ra 0, $$
$$0 \ra R(-a-3)^{a} \ra R(-a-1)^{3a-3}  \oplus R(-a-2)^3
 \ra  R(-a)^{2a}\oplus R(-3)  \ra I_{Z_{o(a)}} \ra 0,$$
$$0 \ra R(-a-3)^{2a} \ra R(-a-2)^{3a+3} \ra R(-a)^{a} \oplus R(-a-1)^3 \oplus R(-3)  \ra I_{Z_{p(a)}}\ra 0.$$
As a consequence, the Minimal Resolution Conjecture on $X$  holds for $m(a)$, $n(a)$, $o(a)$ and $p(a)$ whenever $a \geq 1$.
\end{teo}

Our results are based on the proof of a Proposition of Hartshorne (\cite{Hargor} Proposition 2.4), where it is proved that any set of
general points on $X$ can be $G$-linked to a point in $X$ by a
finite number of $G$-links. The main idea used to prove that
result was the following Lemma.
\begin{lema}\label{lligam} Let $d,g$ be the degree and genus of a smooth arithmetically Cohen-Macaulay curve $D$ on $X$ and let  $n,n' \in \NN$ satisfy $g \leq n,n' \leq d+g-1$ and
$n+n'=deg(mH_D-K_D)$ for some $m$. Then $n$ (respectively $n'$) general points
on $X$ lie on a smooth arithmetically Cohen-Macaulay curve $C$ in the linear system $|D|$ and a there is a divisor $G \subset C$
linearly equivalent to $mH_C-K_C$ which links $n$ general points to
$n'$ general points.
\end{lema}
\begin{proof} See the proof of \cite{Hargor} Proposition 2.4.
\end{proof}

\begin{rk}\rm a) In applying Lemma \ref{lligam} it will be important to keep in mind that, for an arithmetically Cohen-Macaulay curve $C$ on $X$, the couple $(d,g)$ of degree and genus of $C$ determines the linear system in which $C$ belongs (see the paragraph preceding Theorem \ref{MRC}). 
\\
b) In the hypothesis of Lemma \ref{lligam}, $n$
(respectively $n'$) general points on $X$ form a general divisor
of degree $n$ on $C.$ However, the results
of \cite{Mustata} for general points on curves do not apply to prove the Minimal Resolution Conjecture for $m(a)$, $n(a)$, $o(a)$ or $p(a)$.
\end{rk}

We follow the proof of \cite{Hargor} Proposition 2.4 in reverse
order: we perform links to obtain larger sets of points. In this
way and according to \cite{Hargor} 2.4, we can pass from $Z_{m(a-1)}$
general points on $X$ to $Z_{m(a+1)}$ general points by the following
sequence of $G$-links in $X$:

\begin{itemize}\label{paglinks}
\item First link: $Z_{m(a-1)}$ is $G$-linked to $Z_{n(a-1)}$ by
an arithmetically Gorenstein set of points $G$ linearly equivalent
to $(a-1)H_C$ on a smooth curve $C$ linearly equivalent to $(a-1)H$ on
$X$.
\item Second link: $Z_{n(a-1)}$ is $G$-linked to $Z_{o(a)}$ by
an arithmetically Gorenstein set of points $G$ linearly equivalent
to $(2a-2)H_C-K_C$ on a smooth curve $C$ linearly equivalent to
$C_0+(a-1)H$ on $X$.
\item Third link: $Z_{o(a)}$ is $G$-linked to $Z_{p(a)}$ by
an arithmetically Gorenstein set of points $G$ linearly equivalent
to $(2a-1)H_C-K_C$ on a smooth curve $C$ linearly equivalent to
$\Gamma+(a-1)H$ on $X$.
\item Fourth link: $Z_{p(a)}$ is $G$-linked to $Z_{m(a+1)}$ by
an arithmetically Gorenstein set of points $G$ linearly equivalent
to $2aH_C-K_C$ on a smooth curve $C$ linearly equivalent to
$C_0+aH$ on $X$.
\end{itemize}

This sequence of linkages will allow us to deduce the minimal free
resolution of $Z_{m(a+1)}$ from that of $Z_{m(a-1)}$. To do this we need
the following four results corresponding to the four links
mentioned above:

\begin{prop}[First link]\label{firstlink} Assume that the ideal $I_{Z_{m(a)},X} \subset R_X$ of $m(a)=\frac{3}{2}a(a-1)+a$
general points on $X$ has the following
graded minimal free resolution over $R$
$$0 \ra R(-a-3)^{a-1} \ra R(-a-1)^{3a}  \ra R(-a)^{2a+1}  \ra I_{Z_{m(a)},X}
 \ra 0. $$ Then the ideal $I_{Z_{n(a)},X} \subset S$ of $n(a)=\frac{3}{2}a(a-1)+2a$ general points on $X$ has the
 following minimal free resolution over $R$:
$$0 \ra R(-a-3)^{2a-1} \ra R(-a-2)^{3a}  \ra R(-a)^{a+1}  \ra I_{Z_{n(a)}, X}
 \ra 0 . $$
\end{prop}

\begin{prop}[Second link]\label{secondlink}
Assume that the ideal $I_{Z_{n(a-1)},X} \subset R_X$ of
$n(a-1)=\frac{3}{2}(a-2)(a-1)+2a-2$ general points on $X$ has the
 following minimal free resolution over $R$
$$0 \ra R(-a-2)^{2a-3} \ra R(-a-1)^{3a-3}  \ra R(-a+1)^{a}  \ra I_{Z_{n(a-1)}, X}
 \ra 0 . $$ Then the ideal $I_{Z_{o(a)},X}$ of $o(a)=\frac{3}{2}a(a-1)+a+1$ general points on $X$ has the following
 minimal free resolution over $R:$
$$0 \ra R(-a-3)^{a} \ra R(-a-1)^{3a-3}  \oplus R(-a-2)^3
 \ra  R(-a)^{2a} \ra I_{Z_{o(a)},X} \ra 0 .$$
\end{prop}

\begin{prop}[Third link]\label{thirdlink} Assume that the ideal $I_{Z_{o(a)},X}$ of $o(a)=\frac{3}{2}a(a-1)+a+1$ general points on $X$ has the following
 minimal free resolution over $R:$
\begin{multline*}0 \ra R(-a-3)^{a} \ra R(-a-1)^{3a-3}  \oplus R(-a-2)^3
\\
\ra  R(-a)^{2a} \ra I_{Z_{o(a)},X} \ra 0 .\end{multline*} Then the
ideal $I_{Z_{p(a)},X}$ of $p(a)=\frac{3}{2}a(a-1)+2a+1$ general points
on $X$ has a minimal free resolution over $R$:
$$0 \ra R(-a-3)^{2a} \ra R(-a-2)^{3a+3} \ra R(-a)^{a} \oplus R(-a-1)^3  \ra I_{Z_{p(a)},X}\ra 0.$$
\end{prop}

\begin{prop}[Fourth link]\label{fourthlink}
Assume that the ideal $I_{Z_{p(a)},X}$ of $p(a)=\frac{3}{2}a(a-1)+2a+1$
general points on $X$ has the following minimal free
$R$-resolution:
$$0 \ra R(-a-3)^{2a} \ra R(-a-2)^{3a+3} \ra R(-a)^{a} \oplus R(-a-1)^3  \ra I_{Z_{p(a)},X}\ra 0.$$ Then the ideal
$I_{Z_{m(a+1)},X}$ of $Z_{m(a+1)}$ general points on $X$ has a minimal
free resolution over $R$:
$$0 \ra R(-a-4)^{a} \ra R(-a-2)^{3a+3}  \ra R(-a-1)^{2a+3}  \ra I_{Z_{m(a+1)},X}
 \ra 0. $$
\end{prop}

Assuming that we had proven the four propositions we prove the
theorem now:

\begin{proofteo} 
For $a=1,2$ and $3$ we are considering the following number of general points on $X$:
$m(1)=1, n(1)=o(1)=2, p(1)=3, m(2)=5, n(2)=7, o(2)=6,p(2)=8, m(3)=12,n(3)=15,o(3)=13,p(3)=16.$
These are also general points in $\PP^3$ and therefore the Minimal Resolution Conjecture is known to hold. Moreover, for $a=3$ their minimal free resolutions correspond to the ones given in the statement of the Theorem. 

For $a \geq 4$ we shall first prove the following claim: 
 
\emph{Claim:} For $a \geq 2$ the resolution of $Z_{m(a)}$ in $X$ is the one given in Proposition \ref{firstlink}, i.e.
\begin{equation}\label{hipot}
0 \ra R(-a-3)^{a-1} \ra R(-a-1)^{3a}  \ra R(-a)^{2a+1}  \ra I_{Z_{m(a)},X}
 \ra 0.
\end{equation}
\emph{Proof of Claim:}
We proceed by induction on $a$ and we start by proving that it is true for $a=2$ and $a=3$.

$a=2$: $m(2)$ general points in $X$ are 5 general points in
$\PP^3$. Therefore  $Z_{m(2)}$ has the following graded minimal free
resolution:
$$0 \ra R(-5) \ra R(-3)^5 \ra R(-2)^5 \ra I_{Z_{m(2)}} \ra 0.$$ Using the mapping cone procedure applied to the free resolutions
of the ideals in the exact sequence
$$0 \ra I_X \ra I_{Z_{m(2)}} \ra I_{Z_{m(2)},X} \ra 0,$$
we obtain that the minimal free resolution of $I_{Z_{m(2)},X}$ coincides with (\ref{hipot}) for $a=2$.

$a=3$: $m(3)=12$ general points in $\PP^3$ lie on a smooth cubic
surface. Their minimal free resolution is
$$0 \ra R(-6)^2 \ra R(-4)^9 \ra R(-3)^8 \ra I_{Z_{m(3)}} \ra 0. $$
Then we apply the mapping cone construction to the following
commutative diagram:
$$\begin{array}{cccccclll}
0 & \ra & I_X & {\lra} & I_{Z_{m(3)}}  & \ra  & I_{Z_{m(3)},X} & \ra & 0\\
 & & \uparrow & & {\uparrow} & & & & \\
 & & R(-3) & \srel{\xi_1}{\lra} & R(-3)^8  & & & &\\
 & &  {\uparrow} & & {\uparrow} & & &  &\\
 & & 0 & & R(-4)^9  & & & & \\
 & &  & &  {\uparrow} & & & & \\
  & & & & R(-6)^{2} & & & &\\
  & &  & & \uparrow & & & &\\
    & &  & & 0 & & & & \\
   \end{array}$$
As $f$ is one of the generators of $I_{Z_{m(3)}}$, when we write the
map $\xi_1$ as a matrix we see that it contains a non-zero scalar
entry. Therefore, in the resolution of $I_{Z_{m(3)},X}$ obtained by
mapping cone, there is one term $R(-3)$ that can be split off,
giving the desired minimal free resolution.

$a \geq 4$: We assume that the resolution of $Z_{m(a-2)}$ in $X$ is (\ref{hipot}) for $a-2$. 
Then Proposition \ref{firstlink} holds for $a-2$, and so $n(a-2)$ satisfies the hypothesis of Proposition
\ref{secondlink}. This in turn implies that the hypothesis of Proposition
\ref{thirdlink} is satisfied for $a-1$ and that of Proposition \ref{fourthlink}
also holds for $a-1$. Therefore we obtain the desired minimal free
resolution for $I_{Z_{m(a)},X}$ and hence the claim is proved. 

\vspace{2mm}

Now the Theorem immediately follows for $m(a)$ by applying the horseshoe lemma (cf. \cite{Weibel} 2.2.8) to the
following diagram $$\begin{array}{cccccccll}
0 & \ra & I_X & \lra & I_{Z_{m(a)}} & \lra  & I_{Z_{m(a)},X} & \lra & 0\\
 & & \uparrow & &  & & {\uparrow} & & \\
 & & R(-3) & &   & & R(-a)^{2a+1} & &\\
 & & {\uparrow} & &  & & {\uparrow} &  &\\
 & & 0 &  &  & & R(-a-1)^{3a}  & & \\
 & &  & &  & &{\uparrow} & & \\
  & & & &  & & R(-a-3)^{a-1} & &\\
  & &  & & & &\uparrow  & &\\
    & &  & &  & &0 & &   \end{array}$$
The resolution obtained is minimal for $a \geq 3$ (but not for $a=2$).

We now prove the Theorem for $n(a), o(a)$ and $p(a)$. The claim above and Propositions \ref{firstlink}, \ref{secondlink},
\ref{thirdlink} and \ref{fourthlink} provide us minimal free resolutions for $I_{Z_{n(a)},X}$ when $a \geq 2$ and for $I_{Z_{o(a)},X}$, $I_{Z_{p(a)},X}$ when $a \geq 3$. The horseshoe lemma gives us the desired minimal free resolutions for $n(a), o(a)$ and $p(a)$ when $a \geq 3$.
\end{proofteo}

The hard work is proving the four propositions.

\begin{prooflema1}
$n=m(a)$ and $n'=n(a)$ satisfy the hypotheses of Lemma \ref{lligam} with $d=3a$ and $g=1/2(3a^2-3a+2)$ (which are the degree and genus of a curve $C \sim aH$) because $n+n'=deg(aH_C)$. Therefore $Z_{m(a)}$ and $Z_{n(a)}$ are directly
$CI$-linked on a smooth curve $C\sim aH$ by $G \sim aH_C$ and $I_G \subset R$ has the following minimal free resolution $$0 \ra R(-2a-3) \ra R(-2a) \oplus R(-a-3)^2 \ra R(-a)^2 \oplus R(-3) \ra I_G \ra 0.$$
By mapping cone construction applied to the commutative diagram
$$\begin{array}{cccccclll}
0 & \ra &  I_{G} & \srel{i}{\lra} & I_{Z_{m(a)}} & \lra  & I_{Z_{m(a)}}/I_{G}  \ra  0\\
 & & e_1 \uparrow & & d_1{\uparrow} & &  \\
 & & R(-a)^2 \oplus R(-3) & \srel{\xi_1}{\lra} & R(-a)^{2a+1} \oplus R(-3)   & & \\
 & & e_2 {\uparrow} & & d_2 {\uparrow} & &\\
 & & R(-2a) \oplus R(-a-3)^2  & \srel{\xi_2}{\lra} & R(-a-1)^{3a}  & & \\
 & & e_3 \uparrow & & d_3 {\uparrow} & & \\
  & & R(-2a-3) & \srel{\xi_3}{\lra} & R(-a-3)^{a-1} & & \\
  & & \uparrow & & \uparrow & & \\
    & & 0 & & 0 & &  \\
   \end{array}$$
we obtain the following free resolution of $I_{Z_{n(a)}}$ (see Lemma
\ref{mapcone})
\begin{multline*}
0 \ra R(-a-3)^{2a+1} \oplus R(-2a) \ra R(-a-2)^{3a} \oplus
R(-a-3)^2 \oplus R(-2a) \\
\ra R(-a)^{a+1} \oplus R(-3) \ra I_{Z_{n(a)}} \ra 0.
\end{multline*}
 The
generators of $I_G$ can be taken in the set of minimal generators
of $I_{Z_{m(a)}}$, so the map $\xi_1$ has non-zero scalar entries. This
implies that the summand $R(-a-3)^2 \oplus R(-2a)$ in the free
resolution above can be split off. Therefore, a minimal free
resolution for $I_{Z_{n(a)}}$ is
$$0 \ra R(-a-3)^{2a-1} \ra R(-a-2)^{3a} \ra R(-a)^{a+1} \oplus R(-3) \ra I_{Z_{n(a)}}  \ra 0.$$
The mapping cone construction applied now to the exact sequence $0
\ra I_X \ra I_{Z_{n(a)}} \ra I_{Z_{n(a)},X} \ra 0$ gives the desired
minimal free resolution for $I_{Z_{n(a)},X}.$
\end{prooflema1}

Before proving Proposition \ref{secondlink} we need the following
Lemma.
\begin{lema}\label{resGor} Let $C$ be a curve in $X$ linearly equivalent to $C_0+aH$ and let $G$ be an effective divisor on $C$ linearly equivalent to
$2aH_C-K_C$. Then the minimal free $R$-resolution of $I_{G,C}$ is
$$0 \ra R(-2a-4) \ra R(-a-3) \oplus R(-a-2) \oplus R(-2a-1)
\ra R(-a-1)^2 \ra I_{G,C} \ra 0. $$
\end{lema}
\begin{proof} As $C \sim C_0+aH$, we have that  $I_{C,X} \cong
H^0_{\ast}(X,\JJ_{C_0}(-a))$ and therefore we know the $R$-
resolution of $I_{C,X}$. Thus we can apply horseshoe lemma to the
following diagram
$$\begin{array}{cccccccll}
0 & \ra & I_X & \lra & I_C & \lra  & I_{C,X} & \lra & 0\\
 & & \uparrow & &  & & {\uparrow} & & \\
 & & R(-3) & &   & & R(-a-1) \oplus R(-a-2) & &\\
 & & {\uparrow} & &  & & {\uparrow} &  &\\
 & & 0 &  &  & & R(-a-3)^{2}  & & \\
 & &  & &  & &{\uparrow} & & \\
  & &  & &  & & 0 & &\\
   \end{array}$$
  and we obtain the following minimal free $R$-resolution for
  $R/I_C$:
  $$0 \ra R(-a-3)^2 \ra R(-a-1) \oplus R(-a-2) \oplus R(-3) \ra R \ra R/I_C \ra 0.$$
As $C$ is an arithmetically Cohen-Macaulay curve, applying
Hom(.,R(-4)) to this exact sequence we obtain a minimal free
resolution for the canonical module
$\omega_C=Ext^2_R(R/I_C,R(-4))$:
$$0 \ra R(-4) \ra R(a-3) \oplus R(a-2) \oplus R(-1) \ra R(a-1)^2 \ra \omega_C \ra 0. $$
When we shift this resolution by $-2a$ we obtain the desired
minimal free resolution for $H^0_{\ast}(\OO_C(K_C-2aH_C)) \cong
I_{G,C}$:
$$0 \ra \OO_{\PP^3}(-2a-4) \ra R(-a-3) \oplus R(-a-2) \oplus R(-2a-1)
\ra R(-a-1)^2 \ra I_{G,C} \ra 0. $$
\end{proof}

\begin{prooflema2}
$n=n(a-1)$ and $n'=o(a)$ satisfy the hypotheses of Lemma
\ref{lligam} with $d=3a-1$, $g=1/2(3a^2-5a+2)$ because $n+n'= deg( (2a-2)H_C-K_C)$ for a curve $C \sim C_0+(a-1)H$. Therefore there is a smooth curve $C\sim C_0+(a-1)H$ so that $n(a-1)$ general
points on $X$ and $o(a)$ general points on $X$ are $G$-linked by $G \sim (2a-2)H_C-K_C$ on $C$.

According to lemma \ref{resGor}, the minimal free resolution of
$I_{G,C}$ is
$$0 \ra R(-2a-2) \ra R(-a-2) \oplus R(-a-1) \oplus R(-2a+1)
\ra R(-a)^2 \ra I_{G,C} \ra 0. $$

On the other hand, the mapping cone procedure applied to the
following exact sequence
$$0 \ra I_{C,X} \ra I_{Z_{n(a-1)},X} \ra I_{Z_{n(a-1)},C} \ra 0$$
gives us the minimal free $R$-resolution of $I_{Z_{n(a-1)},C}$:
$$0 \ra R(-a-2) \ra R(-a) \oplus R(-a-1)^{3a-2} \ra R(-a+1)^{a} \ra I_{Z_{n(a-1)},C} \ra 0.$$
Now we apply the mapping cone construction to the following
commutative diagram:
$$\begin{array}{cccccclll}
0 & \ra &  I_{G,C} & \srel{i}{\lra} & I_{Z_{n(a-1)},C} & \lra  & I_{Z_{n(a-1)},C}/I_{G,C}  \ra  0\\
 & & e_1 \uparrow & & d_1{\uparrow} & &  \\
 & & R(-a)^2 & \srel{\xi_1}{\lra} & R(-a+1)^{a}   & & \\
 & & e_2 {\uparrow} & & d_2 {\uparrow} & &\\
 & & R(-a-1) \oplus R(-a-2)  & \srel{\xi_2}{\lra} & R(-a-1)^{3a-2} \oplus R(-a) & & \\
 & & \oplus R(-2a+1) & &  & & \\
 & & e_3 \uparrow & & d_3 {\uparrow} & & \\
  & & R(-2a-2) & \srel{\xi_3}{\lra} & R(-a-2)^{2a-1} & & \\
  & & \uparrow & & \uparrow & & \\
    & & 0 & & 0 & &  \\
   \end{array}$$
Writing $\xi_2$ in a matrix way we shall prove that the first
column of $\xi_2$ has a non-zero scalar entry. We assume that the
first column of $\xi_2$ is
$$\begin{pmatrix}
l \\
0\\
\vdots\\
0
\end{pmatrix}$$
for some linear form $l$, and we will get to a contradiction.
Choosing basis, we may denote by $\begin{pmatrix}
L_1 \\
L_2\\
\end{pmatrix}$ the first column of $e_2$, by $\begin{pmatrix}
x_1 \\
x_2\\
\vdots\\
x_a
\end{pmatrix}$
the first column of $d_2$, and the entries of $\xi_1$ by
$$\begin{pmatrix}
p_1 & q_1 \\
p_2 & q_2\\
\vdots & \vdots\\
p_a & q_a
\end{pmatrix}.$$
Then the equality $d_2\xi_2=\xi_1e_2$ implies that
\begin{equation}\label{eqxi}
x_il=p_iL_1+q_iL_2
\end{equation}
for all $i \in \{1, \dots, a\}$. Note that, according to the proof
of Lemma \ref{resGor}, $L_1,L_2$ are given by the curve $C.$

We now consider two different possibilities: either $l,L_1,L_2$
form a complete intersection in $R$ or either $l=\al L_1+\be L_2$
for some $\al, \be \in k$.

If $l,L_1,L_2$ form a complete intersection, then ($\ref{eqxi}$)
implies that $p_i \in (l,L_2)$ and $q_i \in (l,L_1)$ for all $i
\in \{1, \dots, a\}$. Thus, if
$(\overline{f_1},\dots,\overline{f_a}) \subset R/I_C$ are the
generators of $I_{Z_{n(a-1)},C}$, we obtain that both generators of
$I_{G,C}$ belong to $(\overline{l},\overline{L_1},\overline{L_2})
\subset R/I_C$ because they are defined by the product
$$(\overline{f_1},\dots,\overline{f_a})\begin{pmatrix}
p_1 & q_1 \\
p_2 & q_2\\
\vdots\\
p_a & q_a
\end{pmatrix}.$$
But this implies that $I_{G,C}\subseteq
(\overline{l},\overline{L_1},\overline{L_2})$ which is impossible
because $G$ is formed by two sets of general points in $C$.

If $l=\al L_1+\be L_2$, then we can change the generators of
$I_{G,C}=(\overline{x},\overline{y})$ in order to make $L_1$ or
$L_2$ equal to $l$. Let us assume $L_2=l$ and the other case can
be treated analogously. From equation (\ref{eqxi}) we obtain in
this case that $p_i \in (L_2)$ for all $i \in \{1, \dots, a\}$ and
in particular $I_{G,C} \subset (\overline{L_2},\overline{y})$. As
$\overline{L_2},\overline{y}$ do not form a complete intersection
in $R/I_C$, there is an associated prime of
$(\overline{L_2},\overline{y})$ of height 1. But this associated
prime would contain a prime appearing in the primary decomposition
of $I_{G,C}$ and this is impossible because all of them have
height 1 in $R/I_C$.

In both cases we obtain a contradiction. Therefore there is a
non-zero scalar entry in $\xi_2$ and, in particular, in the
resolution of $I_{Z_{n(a-1)},C}/I_{G,C} \cong I_{Z_{n(a-1)}}/I_G$
obtained by mapping cone, one term $R(-a-1)$ can be split off.
Lemma \ref{mapcone} gives then the desired minimal free resolution
for $I_{Z_{o(a)}}.$
\end{prooflema2}

Before proving proposition \ref{thirdlink} we need the following
result:
\begin{lema}\label{lemacorba} Assume that the ideal $I_{Z_{o(a)},X}$ of $o(a)=\frac{3}{2}a(a-1)+a+1$ general points $Z_{o(a)}$ on $X$ has the following
 minimal free resolution over $R:$
\begin{multline*}0 \ra R(-a-3)^{a} \ra R(-a-1)^{3a-3}  \oplus R(-a-2)^3
\\
\ra  R(-a)^{2a} \ra I_{Z_{o(a)},X} \ra 0 .\end{multline*} Then $Z_{o(a)}$
lies on a smooth curve $C \sim \Gamma+(a-1)H$ on $X$ and
$H^0_{\ast}(\JJ_{Z_{o(a)},C})$ has the following minimal free
resolution:
$$ 0 \ra R(-a-3)^{a} \ra R(-a-1)^{3a} \ra R(-a)^{2a} \ra H^0_{\ast}(\JJ_{Z_{o(a)},C}) \ra 0$$
\end{lema}

\begin{proof} 
As $dim|\Gamma+(a-1)H|=3a+1/2(3a^2-3a)-1 \geq o(a)$, any set of $o(a)$ general points lies on a smooth curve $C \sim \Gamma +(a-1)H$.
As $C$ is an arithmetically Cohen-Macaulay curve, we have the following exact sequence
$$0 \ra H^0_{\ast}(\JJ_{C,X}) \srel{i}{\ra} H^0_{\ast}(\JJ_{Z_{o(a)},X}) \ra H^0_{\ast}(\JJ_{Z_{o(a)},C}) \ra 0.$$ We can lift the map
$i$ to the minimal free resolutions of $I_{Z_{o(a)},X}$ and $I_{C,X}$
to get a commutative diagram:

$$\begin{array}{cccccclll}
0 & \ra & H^0_{\ast}(\JJ_{C,X}) & \srel{i}{\lra} & H^0_{\ast}(\JJ_{Z_{o(a)},X}) & \ra  & H^0_{\ast}(\JJ_{Z_{o(a)},C}) & \ra & 0\\
 & & \uparrow & & d_1{\uparrow} & & & & \\
 & & R(-a-1)^3 & \srel{\xi_1}{\lra} & R(-a)^{2a}  & & & &\\
 & & \Phi {\uparrow} & & d_2 {\uparrow} & & &  &\\
 & & R(-a-2)^3 & \srel{\xi_2}{\lra} & R(-a-2)^3 \oplus R(-a-1)^{3a-3} & & & & \\
 & & \uparrow & & d_3 {\uparrow} & & & & \\
  & & 0 & & R(-a-3)^{a} & & & &\\
  & &  & & \uparrow & & & &\\
    & &  & & 0 & & & & \\
   \end{array}$$
and the mapping cone construction gives a free resolution  $$ 0
\ra R(-a-3)^{a} \oplus R(-a-2)^3 \ra R(-a-1)^{3a} \oplus R(-a-2)^3
\ra R(-a)^{2a} \ra H^0_{\ast}(\JJ_{Z_{o(a)},C}) \ra 0 \, ,$$ not
necessarily minimal. In order to prove that the terms $R(-a-2)^3$
can be split off we need to check that, if after choosing basis we
write $\xi_2 =
  \begin{pmatrix}
    \sigma  \\
   \Lambda
  \end{pmatrix}
$  where $\sigma$ is a $3 \times 3$ matrix of scalars and
$\Lambda$ is a $(3a-3) \times 3$ matrix of linear homogeneous
entries, then $\sigma$ has rank 3.

Let us assume that $\sigma$ has rank $\leq 2$ instead. Then there
is a combination of its columns that gives the 0 column. We can
assume (after changing basis if necessary) that the first column
of $\sigma$ is zero. We write the matrix corresponding to $d_2$ as
$
  \begin{pmatrix}
    \delta & D
  \end{pmatrix}
$ (where $\delta$ is a $2a \times 3$ matrix and $D$ a $2a\times
(3a-3)$ matrix), the first column of $\Lambda$ as $\begin{pmatrix}
    l_1 \\
    l_2 \\
    \vdots \\
    l_{3a-3}
  \end{pmatrix}$ and the first column of $\Phi$ as
$\begin{pmatrix} \varphi_1 \\ \varphi_2\\ \varphi_3
\end{pmatrix}$.

Then from the equality $$
  \begin{pmatrix}
    \delta & D
  \end{pmatrix}
  \begin{pmatrix}
    \sigma  \\
    \Lambda
  \end{pmatrix}
= \xi_1 \Phi
$$
applied to the first column of $\xi_2$ and $\Phi$ we obtain that
\begin{equation}\label{ig0}
\overline{D}\cdot
  \begin{pmatrix}
   \overline{ l_1 }\\
   \overline{ l_2} \\
    \vdots \\
   \overline{ l_{3a-3}}
  \end{pmatrix}=0 \textrm{ in }  R/ (\varphi_1,\varphi_2,\varphi_3).
\end{equation}
As the matrix $\Phi$ has determinant $f$, the ideal
$(\varphi_1,\varphi_2,\varphi_3)$ might have either height 2 or 3
in $R$. We consider each case separately:

\underline{Case 1:} height $(\varphi_1,\varphi_2,\varphi_3)=3.$

In this case, in  the ring $R/ (\varphi_1,\varphi_2,\varphi_3)$
each $\overline{ l_i}$ is equal to some linear form $\bar l$
multiplied by some constant $c_{i}$. Therefore (\ref{ig0}) becomes
$$\overline{l} \overline{D}\cdot
  \begin{pmatrix}
    c_1 \\
    c_2 \\
    \vdots \\
    c_{3a-3}
  \end{pmatrix} =0 \textrm{ in } R/
(\varphi_1,\varphi_2,\varphi_3).$$ As $\overline{l}$ is a non-zero
divisor in $R/ (\varphi_1,\varphi_2,\varphi_3)$,we have that
$$
\overline{D}\cdot
  \begin{pmatrix}
    c_1 \\
    c_2 \\
    \vdots \\
    c_{3a-3}
  \end{pmatrix} =0 \textrm{ in } R/
(\varphi_1,\varphi_2,\varphi_3).$$
 This means that there exists a
$k$-linear combination of the columns of $\bar D$ that gives the 0
column in $R/ (\varphi_1,\varphi_2,\varphi_3).$ Changing basis if
necessary, we can assume that there is a column in the matrix $D$
whose entries  belong to the ideal
$(\varphi_1,\varphi_2,\varphi_3).$

Let $p$ be the point corresponding to the ideal
$\mathfrak{p}=(\varphi_1,\varphi_2,\varphi_3)$ (note that this is
a point in $X$ because $f$ is the determinant of $\Phi$, so $f \in
\mathfrak{p}$). We split the resolution of $I_{Z_{o(a)},X}$ into short
exact sequences:

\begin{equation}\label{seq1} 0 \ra R(-a-3)^{a} \ra R(-a-1)^{3a-3}  \oplus R(-a-2)^3
\srel{d_2}{\ra} E \ra 0 \end{equation}
\begin{equation}\label{seq2} 0 \ra E \ra  R(-a)^{2a} \ra I_{Z_{o(a)},X} \ra 0 \end{equation}
and we are going to check how these sequences behave when we
tensor by $R_\mathfrak{p}$. First of all we need to prove the
following claim, though:

\underline{Claim:} $p$ is not a point in $Z_{o(a)}$

\underline{Proof of claim:} If $p$ was a point in $Z_{o(a)}$ then it
would be a point in the curve $C$. Then, tensoring the minimal
free resolution of $H^0_{\ast}(\JJ_{C,X})$ by $R_\mathfrak{p}$, we
obtain that $\Phi_\mathfrak{p}$ is an isomorphism. But this is a
contradiction since the entries of a column of $\Phi$ belong to
$\mathfrak{p}$, so the claim is proved.

\vspace{2mm}

As $p \notin Z_{o(a)}$, we have that $R_X/I_{Z_{o(a)},X}
{\otimes}_{R}R_\mathfrak{p}=0$ and in particular
$$Tor^R_i(R_X/I_{Z_{o(a)},X},R_\mathfrak{p})=0$$ 
for $i \geq 0.$ From
this we also obtain that $I_{Z_{o(a)},X} {\otimes}_{R}R_\mathfrak{p}$
is isomorphic to $R_X {\otimes}_{R}R_\mathfrak{p}.$

From the exact sequence
$$0 \ra I_X=(f) \ra R \ra R_X \ra 0$$
we obtain that $Tor^R_1(R_X,R_\mathfrak{p})=0.$

Now localizing the exact sequence
$$0 \ra I_{Z_{o(a)},X} \ra R_X \ra R_X/I_{Z_{o(a)},X} \ra 0 $$ and using the
vanishing of the Tor groups above,  we get that
$Tor_i(I_{Z_{o(a)},X},R_\mathfrak{p})=0$ for $i=1,2.$

Thus, after tensoring sequence (\ref{seq2}) by
$\otimes_{R}R_\mathfrak{p}$ we obtain the following short exact
sequence:
\begin{equation}\label{seq3} 0 \ra
E{\otimes}_{R}R_\mathfrak{p} \ra R_\mathfrak{p}^{2a} \ra
I_{Z_{o(a)},X}{\otimes}_{R}R_\mathfrak{p} \ra 0
\end{equation}

Moreover, since $Tor_2(I_{Z_{o(a)},X},R_\mathfrak{p})=0$ and
$Tor^R_1(R,R_\mathfrak{p})=0$, we obtain from sequence
(\ref{seq2}) that $Tor_1(E,R_\mathfrak{p})=0.$ Therefore, we have
the following short exact sequence after localizing sequence
(\ref{seq1}):
\begin{equation}\label{seq4}
0 \ra R_\mathfrak{p}^{a} \ra R_\mathfrak{p}^{3a}
\srel{d_2{\otimes}_{R}R_\mathfrak{p}}{\lra}
E{\otimes}_{R}R_\mathfrak{p}  \ra 0.
\end{equation}
Now we are going to prove that $E{\otimes}_{R}R_\mathfrak{p} $ is
a free $R_\mathfrak{p}$-module. To prove this we need to check
that the local cohomology modules
$H^i_{\mathfrak{m}_\mathfrak{p}}(E{\otimes}_{R}R_\mathfrak{p})$
vanish for $i \leq 2$, where $\mathfrak{m}_\mathfrak{p}$ is the
maximal ideal of the ring $R_\mathfrak{p}$. It follows immediately
from sequence (\ref{seq4}) that these modules are zero for
$i=0,1$. Moreover, from the exact sequence (\ref{seq3}), taking
into account the isomorphism $I_{Z_{o(a)},X}
{\otimes}_{R}R_\mathfrak{p} \cong R_X
{\otimes}_{R}R_\mathfrak{p}$, we obtain that
$H^2_{\mathfrak{m}_\mathfrak{p}}(E{\otimes}_{R}R_\mathfrak{p})=0$
because $depth(R_X {\otimes}_{R}R_\mathfrak{p})=2$.

Therefore, $E{\otimes}_{R}R_\mathfrak{p}$ is a free module and the
exact sequence (\ref{seq4}) splits. But this contradicts the fact
that a column of $d_2$ has all its entries in the ideal
$\mathfrak{p}.$

\underline{Case 2:} height $(\varphi_1,\varphi_2,\varphi_3)=2.$

Let us assume that
$(\varphi_1,\varphi_2,\varphi_3)=(\varphi_1,\varphi_2)$ (the other
cases are similar). Then $(\varphi_1,\varphi_2)$ is the ideal of a
line in $X.$

From (\ref{ig0}) we see that in $R/(\varphi_1,\varphi_2)$
$$
  \begin{pmatrix}
    \overline{l_1} &
    \overline{l_2} &
    \dots &
    \overline{l_{3a-3}}
  \end{pmatrix}
  \overline{D^t} =0
$$
The ideal $(\overline{l_1},\overline{l_2},
\dots,\overline{l_{3a-3}})$ has height at most 2 in
$R/(\varphi_1,\varphi_2)$. Performing a base change if necessary
we can assume that $(\overline{l_1},\overline{l_2},
\dots,\overline{l_{3a-3}})=(\overline{l_1},\overline{l_2})$. The
case where $\overline{l_1} \in (\overline{l_2})$ or viceversa is
easier and left to the reader.

The columns of $\overline{D^t}$ are syzygies of
$(\overline{l_1},\overline{l_2}, \dots,\overline{l_{3a-3}})$. This
means that $\overline{D^t}=Z.T$ where $Z$ is the following
$(3a-3)\times (3a-4)$ matrix
$$
\begin{pmatrix}
\overline{l_2} & z_1^2 & z_1^3 & \dots & z_1^{3a-4} \\
-\overline{l_1} & z_2^2 & z_2^3 & \dots & z_2^{3a-4} \\
0 & z_3^2 & 0 & \dots & 0 \\
0 & 0 &  z_4^3 & \dots & 0 \\
  &   &   \vdots &  & \\
0 & 0 &  0 & \dots & z_{3a-3}^{3a-4}
\end{pmatrix},
$$
$z^i_j \in k$, and $T$ is a $(3a-4) \times (2a)$ matrix of
polynomials in $R/(\varphi_1,\varphi_2)$. Considering the
following constant vector
$$\mu = (1,0,-\frac{z_1^2}{z_3^2},-\frac{z_1^3}{z_4^3}, \dots,
-\frac{z_1^{3a-4}}{z_{3a-3}^{3a-4}}),$$ we have that
$\overline{D}\mu^t=T^t Z^t \mu^t$ is equal to
$$T^t \begin{pmatrix}
-\overline{l_2} \\
0 \\
\vdots \\
0
\end{pmatrix}.$$
Therefore $D \mu^t $ has each entry in the ideal
$(\varphi_1,\varphi_2,l_2)$. After doing a base change we can
assume that the entries of one of the columns in $D$ lie all in
the ideal $(\varphi_1,\varphi_2,l_2)$ (or any other linear form
instead of $l_2$ if the reader bothers about too many base
changes). We consider the point $p$ defined by this ideal which is
a point in the line given by $\varphi_1,\varphi_2$. It cannot be a
point of $Z_{o(a)}$ because the points of $Z_{o(a)}$ were general points in
$X$, so none of them can lie in one of the 27 lines of $X$. The
rest of the proof in this case follows by localizing at $p$,
exactly as we did in case 1.

Therefore, in both cases we obtain that the rank of $\sigma$ must
be 3, which in turn implies that the term $R(-a-2)^3$ in the
resolution of $H^0_{\ast}(\JJ_{Z_{o(a)},C})$ is redundant.
\end{proof}

Now we are ready to prove Proposition \ref{thirdlink}.

\begin{prooflema3}
$n=o(a)$ and $n'=p(a)$ satisfy the hypotheses of Lemma \ref{lligam} with $d=3a$ and $g=1/2(3a^2-3a)$ because $n+n'=deg((2a-1)H_C-K_C)$ for a curve 
$C \sim \Gamma +(a-1)H$. Therefore, there exists a smooth curve  $C \sim \Gamma +(a-1)H$ and an arithmetically Gorenstein scheme $G \sim (2a-1)H_C-K_C$ so that $Z_{o(a)}$ and $Z_{p(a)}$ are linked by $G$ on $C$. Then on the curve $C$, $Z_{o(a)}+Z_{p(a)}$ is linearly
equivalent to $(2a-1)H_C - K_C$. From Lemma \ref{lemacorba} we
have a minimal free resolution of the $R$-module
$M=H^0_{\ast}(\OO_C(-Z_{o(a)}))$:
$$ 0 \ra R(-a-3)^{a} \ra R(-a-1)^{3a} \ra R(-a)^{2a} \ra H^0_{\ast}(\OO_C(-Z_{o(a)})) \ra 0.$$
As $M$ is a Cohen-Macaulay module we can dualize applying
$\Hom(\cdot,R(-4))$ and we obtain a minimal free resolution of
$Ext^2(M,R(-4))$: \begin{equation}\label{resdual} 0 \ra
R(a-4)^{2a} \ra R(a-3)^{3a} \ra R(a-1)^a \ra Ext^2(M,R(-4)) \ra 0.
\end{equation} By local duality $Ext^2(M,R(-4))$ is isomorphic to
$H^2_{\mathfrak{m}}(M)^{\vee}$, which at the same time is
isomorphic to $H^1_{\ast}(\OO_C(-Z_{o(a)}))^{\vee}.$ By Serre duality
we have $H^1_{\ast}(\OO_C(-Z_{o(a)}))^{\vee} \cong
H^0_{\ast}(\OO_C(Z_{o(a)}+K_C))$. But as $Z_{o(a)} +K_C \sim
-Z_{p(a)}+(2a-1)H_C$, twisting (\ref{resdual}) by $R(-2a+1)$ we obtain
a minimal free resolution of $H^0_{\ast}(\OO_C(-Z_{p(a)}))$. Now we
apply the horseshoe lemma to the following diagram:

$$\begin{array}{cccccccll}
0 & \ra & H^0_{\ast}(\JJ_{C,X}) & \lra & H^0_{\ast}(\JJ_{Z_{p(a)},X}) & \lra  & H^0_{\ast}(\JJ_{Z_{p(a)},C}) & \lra & 0\\
 & & \uparrow & &  & & {\uparrow} & & \\
 & & R(-a-1)^3 & &   & & R(-a)^{a} & &\\
 & & {\uparrow} & &  & & {\uparrow} &  &\\
 & & R(-a-2)^3 &  &  & & R(-a-2)^{3a}  & & \\
 & & \uparrow & &  & &{\uparrow} & & \\
  & & 0 & &  & & R(-a-3)^{2a} & &\\
  & &  & & & &\uparrow  & &\\
    & &  & &  & &0 & & \\
   \end{array}$$
and we obtain the desired minimal free resolution (note that it is
already minimal because no terms can be split off.)
\end{prooflema3}

Before proving Proposition \ref{fourthlink} we need the following
result.

\begin{lema}\label{resenC} Let $C \subset X$ be a smooth curve linearly equivalent to
$C_0+aH$ containing $Z_{p(a)}$. If the minimal free resolution of
$I_{Z_{p(a)},X}$ is
$$0 \ra R(-a-3)^{2a} \ra R(-a-2)^{3a+3} \ra R(-a)^{a} \oplus R(a-1)^3  \ra I_{Z_{p(a)},X}\ra 0,$$
then the minimal free resolution of $I_{Z_{p(a)},C}$ is:
$$0 \ra R(-a-3)^{2a+2} \ra R(-a-2)^{3a+4} \ra R(-a)^a \oplus R(-a-1)^2 \ra I_{Z_{p(a)}, C} \ra 0 .$$
\end{lema}
\begin{proof}
We know both the minimal free resolution of $I_{C,X}$ and
$I_{Z_{p(a)},X}$, so we apply the mapping cone procedure to the
following commutative diagram:
$$\begin{array}{cccccclll}
0 & \ra & I_{C,X} & \srel{i}{\lra} & I_{Z_{p(a)},X} & \ra  & I_{Z_{p(a)},C} & \ra & 0\\
 & & e_1 \uparrow & & d_1{\uparrow} & & & & \\
 & & R(-a-1)  \oplus R(-a-2) & \srel{\xi_1}{\lra} & R(-a)^{a} \oplus R(-a-1)^3  & & & &\\
 & & e_2 {\uparrow} & & d_2 {\uparrow} & & &  &\\
 & &  R(-a-3)^2  & \srel{\xi_2}{\lra} & R(-a-2)^{3a+3} & & & & \\
 & & \uparrow & & d_3 {\uparrow} & & & & \\
  & &  0&  & R(-a-3)^{a} & & & &\\
  & & & & \uparrow & & & &\\
    & &  & & 0 & & & & \\
   \end{array}$$
Note that the top row is exact because $H^1_{\ast}(\JJ_{C,X})=0.$
We need to prove that there exists a non-zero scalar entry in the
first column of $\xi_1$. Let us assume that the scalar entries in
$\xi_1$ are zero. Then if
$I_{Z_{p(a)},X}=(\overline{f_1},\dots,\overline{f_a},\overline{f_{a+1}},\overline{f_{a+2}},\overline{f_{a+3}})
\subset R_X $, where $f_{a+1},f_{a+2},f_{a+3}$ are forms of degree
$a+1$, we have that $I_{C,X}$ is generated by two forms
$\overline{x}$, $\overline{y}$ whith $\overline{x} \in
J:=(\overline{f_1},\dots,\overline{f_a}).$ Now if
 $$ \begin{pmatrix}
    Q_1 & Q_2 \\
    L_1 & L_2
  \end{pmatrix}$$
is the matrix corresponding to $e_2$ in certain basis, we have
that $\overline{x}Q_1+\overline{y}L_1$ is 0 in $R/(f)$. Hence, as
$\overline{x}$ is in $J$ we obtain that $\overline{y}L_1 $ is in
$J$ too and thus $(\overline{x},\overline{y}L_1) \subset J$. As
$C$ is an irreducible curve, $\overline{y}$ and $\overline{L_1}$
are coprimes in $R_X=R/(f)$ and thus
$(\overline{x},\overline{y}L_1)=(\overline{x},\overline{y}) \cap
(\overline{x},\overline{L_1}) \subset J $.  But $J$ is formed by
the generators of degree $a$ of the ideal of a \emph{general} set
of points in $C$, so we can deduce that
$(\overline{x},\overline{y}) \subset J$. If we denote by $Z$  the
zero set of $J$, this means that $Z$ is inside $C$. But this
contradicts the fact that $J$ is formed by removing three
generators of the ideal of a general set of points in $C$, and
thus $Z$ cannot lie inside $C$. Therefore there is a non-zero
scalar entry in $\xi_1$, which implies that a summand $R(-a-1)$ is
redundant in the resolution of $I_{Z_{p(a)},C}$ obtained by mapping
cone.
\end{proof}

\begin{prooflema4}
$n=p(a)$ and $n'=m(a+1)$ satisfy the hypotheses of Lemma \ref{lligam}
with $d=3a+2$, $g=\frac{1}{2}(3a^2+a-2)$ (which are the degree and genus of $C \sim C_0+aH$) because $n+n'=deg(2aH_C-K_C)$. Therefore we perform a link on a smooth curve $C \sim C_0+aH$ with an arithmetically Gorenstein scheme $G \sim 2aH_C-K_C$. By Lemma \ref{resGor}, the minimal free resolution of
$I_{G,C}$ is
$$0 \ra R(-2a-4) \ra R(-a-2) \oplus R(-a-3) \oplus R(-2a-1)
\ra R(-a-1)^2 \ra I_{G,C} \ra 0. $$ On the other hand, by Lemma
\ref{resenC} we know the minimal free resolution of $I_{{Z_{p(a)}},C}$,
so we can apply the mapping cone procedure to the following
commutative diagram:
$$\begin{array}{cccccclll}
0 \ra & I_{G,C} & \srel{i}{\lra} & I_{Z_{p(a)},C} & \ra  I_{Z_{p(a)},C}/I_{G,C}   \ra  0\\
  & e_1 \uparrow & & d_1{\uparrow} &  \\
  & R(-a-1)^2 & \srel{\xi_1}{\lra} & R(-a)^{a} \oplus R(-a-1)^2  & \\
  & e_2 {\uparrow} & & d_2 {\uparrow} & \\
  & R(-a-2) \oplus R(-a-3)  & \srel{\xi_2}{\lra} & R(-a-2)^{3a+4} &  \\
  & \oplus R(-2a-1) & &  &  \\
  & e_3 \uparrow & & d_3 {\uparrow} &  \\
  & R(-2a-4) & \srel{\xi_3}{\lra} & R(-a-3)^{2a-2} &\\
  & \uparrow & & \uparrow & \\
  & 0 & & 0 & \\
   \end{array}$$
We shall prove that the first column of $\xi_2$ has a non-zero
scalar entry and that, if we write
$$\xi_1=\begin{pmatrix}
M \\
\sigma
\end{pmatrix}
$$
where $\sigma$ is a $2 \times 2$ matrix with scalar entries, then
$\sigma$ has rank 2.

We start by proving that the first column of $\xi_2$ has a
non-zero scalar entry.

We write $e_2=\left(\begin{array}{ccc} L_1 & Q_1 & H_1\\
L_2 & Q_2 & H_2
\end{array}\right).$
The first two columns form the transpose of the resolution of
$I_{C,X}$ over $R$ (see the proof of Lemma \ref{resGor}), so $L_1$
and $L_2$ are two linear forms defining a line in $X.$
If the first column of $\xi_2$ was 0, then $\xi_1 {\cdot}\left(\begin{array}{c}  L_1\\
L_2
\end{array}\right)=0$ and in particular $\sigma{\cdot}\left(\begin{array}{c}  L_1\\
L_2
\end{array}\right)=0. $ But $L_1$ and $L_2$ are linearly independent over
$k$, so this implies that $\sigma$ is the zero matrix. If we
denote by
$(\overline{f_1},\dots,\overline{f_a},\overline{f_{a+1}},\overline{f_{a+2}})$
the generators of $I_{Z_{p(a)},C}$ in $R/I_C$, where
$deg\overline{f_i}=a$ if $i \leq a$ and
$deg\overline{f_{a+1}}=deg\overline{f_{a+2}}=a+1$, then $\sigma
=0$ implies that $I_{G,C} \subseteq
(\overline{f_1},\dots,\overline{f_a}) \subset I_{Z_{p(a)},C}$. But this
is impossible because $Z_{p(a)}$ is a set of general points in $C$ and
so it cannot be contained in the set of zeroes
$Z(\overline{f_1},\dots,\overline{f_a}) \subset G \subset C.$
Therefore the first column of $\xi_2$ must have a non-zero scalar
entry.

Now we prove that $\sigma$ has rank 2. Indeed, if $\sigma=0$ , we
have just seen how to get to a contradiction. If the rank of
$\sigma$ was 1, then we can assume that the second row of $\sigma$
is 0 (changing basis if necessary). But this implies $I_{G,C}
\subseteq (\overline{f_1},\dots,\overline{f_a},\overline{f_{a+1}})
\subset I_{Z_{p(a)},C}$ and this leads to a contradiction as above.

Therefore, the mapping cone construction tells us that the minimal
free resolution of $I_{Z_{p(a)},C}/I_{G,C} \cong I_{Z_{p(a)}}/I_G$ is
\begin{multline*}
0 \ra R(-2a-4) \ra R(-a-3)^{2a+3} \oplus R(-2a-1) \ra R(-a-2)^{3a+3} \\
\ra R(-a)^{a}\ra I_{Z_{p(a)}}/I_G \ra 0.
\end{multline*}
By Lemma \ref{mapcone}, the resolution of the residual $I_{Z_{m(a+1)}}$ is
obtained by dualizing this resolution and twisting by $R(-2a-4)$.
So this gives the desired free resolution which is minimal because
no terms can be split off.
\end{prooflema4}

\end{document}